\documentclass[aop,MSNbibl,citesort,dvips]{arximspdf}

\doi{10.1214/11-AOP698}
\volume{41}
\issue{2}
\pubyear{2013}
\firstpage{619}
\lastpage{635}

\makeatletter

\newtheorem{prop}{Proposition}
\newtheorem{theorem}{Theorem}
\newtheorem{lem}{Lemma}
\newtheorem{coro}{Corollary}

\newcommand{\T}{\mathbf{T}}

\newcommand{\A}{\mathcal{A}}

\newcommand{\PP}{\mathbb{P}}
\newcommand{\RR}{\mathbb{R}}
\newcommand{\ZZ}{\mathbb{Z}}

\newcommand{\eps}{\varepsilon}
\newcommand{\tiltau}{{\T}}

\makeatother

\begin{document}
\begin{frontmatter}

\title{Recurrence rates and hitting-time distributions
for random walks on the line\thanksref{T1}}
\runtitle{Recurrence rates and hitting-time distributions}

\thankstext{T1}{Supported by the ANR Project TEMI (Th\'eorie ergodique en mesure infinie).}

\begin{aug}
\author[A]{\fnms{Fran\c{c}oise} \snm{P\`ene}\corref{}\ead[label=e1]{francoise.pene@univ-brest.fr}},
\author[A]{\fnms{Beno\^{i}t} \snm{Saussol}\ead[label=e2]{benoit.saussol@univ-brest.fr}}
\and
\author[B]{\fnms{Roland} \snm{Zweim\"uller}\ead[label=e3]{rzweimue@member.ams.org}}
\runauthor{F. P\`ene, B. Saussol and R. Zweim\"uller}
\affiliation{Universit\'e Europ\'eenne de Bretagne and Universit\'e de
Brest, Universit\'e Europ\'eenne de Bretagne and Universit\'e de
Brest, and Universit\"at Wien}
\address[A]{F. P\`ene\\
B. Saussol\\
Laboratoire de Math\'ematiques de Brest\\
Universit\'e Europ\'eenne de Bretagne\\
Universit\'e de Brest\\
CNRS UMR 6205\\
6, Avenue Victor Le Gorgeu\\
CS 93837, 29238 BREST Cedex 3\\
France\\
\printead{e1}\\
\hphantom{E-mail: }\printead*{e2}}
\address[B]{R. Zweim\"uller\\
Fakult\"at f\"ur Mathematik\\
Universit\"at Wien\\
Nordbergstrasse 15, 1090 Wien\\
Austria\\
\printead{e3}} %adresu isvedimo komanda gale!
\end{aug}

% HISTORY:
\received{\smonth{3} \syear{2010}}
\revised{\smonth{7} \syear{2011}}

% ABSTRACT
%
\begin{abstract}
We consider random walks on the line given by a sequence of independent
identically distributed jumps belonging to the strict domain of
attraction of a stable distribution, and first determine the almost
sure exponential divergence rate, as $\varepsilon\to0$, of the return
time to $(-\varepsilon, \varepsilon)$. We then refine this result by
establishing a limit theorem for the hitting-time distributions of
$(x-\varepsilon,x+\varepsilon)$ with arbitrary $x\in\RR$.
\end{abstract}

% KEYWORDS
%

\end{frontmatter}

%s1 ###
\section{Introduction and results}
We consider a recurrent random walk on $\RR$, $S_0:=0$ and
$S_n:=X_1+\cdots+X_n$, $n\geq1$, where
the $X_i$ are i.i.d. random variables on $(\Omega,\mathcal
{F},{\mathbb P})$
such that $\frac{S_n}{A_n}$ converges,
for positive real numbers $A_n$, in distribution
to a stable random variable $X$ with index $\alpha$.
Necessarily (due to recurrence), $\alpha\in[1,2]$, and
the sequence $(A_n)_{n\geq1}$ is regularly varying
of index $\frac{1}{\alpha}$, satisfying
$\sum_{n\geq1}\frac{1}{A_n}=\infty$.\vspace*{1pt}

To capture the speed at which recurrence appears, it is possible to specify,
for such a walk, some deterministic sequences $(\eps_n)$ such that
$S_n\in(-\eps_n,\eps_n)$ infinitely often, or $S_n\notin(-\eps
_n,\eps_n)$
eventually, almost surely. This classical question was addressed, for
example, in~\cite{ChungErdos} and~\cite{BDC}, the results of which
have recently
been extended in~\cite{Cheliotis}.

Here, we are going to study the number of steps it takes
to return to some small neighborhood of the origin
(or to hit a different small interval for the first time).
For related work on random walks in the plane,
intimately related to the $\alpha=1$ case of the
present paper, we refer to~\cite{bf1}.

As an additional standing assumption on our walk, we will always
require the
distribution of the jumps $X_i$ to satisfy the Cram\'er condition
%
%e1 ###
\begin{equation}\label{eqcramer}
{\limsup_{|t|\to\infty}}|{\mathbb E}[e^{itX_1}]|<1.
\end{equation}
This readily implies, in particular, that the event
$\Omega^*:=\{S_n \neq0$ $\forall n\ge1\}$
has positive probability, and $\Omega^*$ has probability one
if and only if no individual path returning to the origin
has positive probability.

As a warm-up we first determine the a.s. rate at which the variables
\[
\tiltau_\varepsilon:=\min\{n\ge1 \dvtx
\vert S_n\vert<\varepsilon\},\qquad \eps>0,
\]
diverge on $\Omega^*$ as $\eps\to0$. Let
$\beta\in[2,\infty]$ be the exponent conjugate to $\alpha$,
that is, $\alpha^{-1}+\beta^{-1}=1$.
\begin{theorem}\label{thmRWP2}
In the present setup,
%
%e2 ###
\begin{equation}\label{eqPointwiseBit}
\lim_{\eps\to0} \frac{\log\tiltau_\varepsilon}{\log\eps}
= -\beta \qquad\mbox{a.s. on } \Omega^{*}.
\end{equation}
\end{theorem}

Our main objective then is to determine the precise order of magnitude
and to
study the asymptotic distributional behavior,
as $\eps\to0$, of the more general hitting times of $\eps$-neighborhoods
of arbitrary given points $x$ on the line. We shall, in fact, do so for
the walk
$S'_n:=S'_0+S_n$, $n\ge0$, with random initial position $S'_0$, independent
of $(S_n)_{n\ge0}$ and having an arbitrary fixed distribution $P$ on
$\RR$.
For any $x\in\RR$ we thus let
\[
\tiltau_\varepsilon^x:=\inf\{m\ge1\dvtx \vert S'_m-x\vert
<\varepsilon\}
\]
and $\Omega^*_x:=\{S'_n\ne x$ $\forall n\ge1\}$.
Outside $\Omega^*_x$ we clearly have $\lim_{\varepsilon\rightarrow0}
\tiltau_\varepsilon^x=\min\{m\ge1 \dvtx S'_m=x\}$.

It is convenient to state the results in terms of, and work
with, the strictly increasing continuous function
$G\dvtx[0,+\infty)\rightarrow[0,+\infty)$ with
$G(0)=0$ which affinely interpolates the values
$G(n)=\sum_{k=1}^n\frac{1}{A_k}$, $n\ge1$.
We denote by $G^{-1}$ its inverse function.
Evidently, $G(n)=o(n)$. Moreover,
by the direct half of Karamata's theorem
(cf. Propositions 1.5.8 and 1.5.9a of~\cite{BGT}),
$G$ is regularly varying with index $\frac{1}{\beta}$,
and satisfies
%
%e3 ###
\begin{equation}\label{eqGasymptotics}\qquad
\frac{n}{A_n}=o(G(n))\qquad
\mbox{if } \alpha=1\qquad
\mbox{while }
\frac{n}{A_n} \sim\frac{G(n)}{\beta}
\mbox{ in case } \alpha\in(1,2].
\end{equation}

We establish a result on convergence in
distribution for $\varepsilon
G(\tiltau_\varepsilon^x)$ conditioned on $\Omega^*_x$
[while $\eps G(\tiltau_\varepsilon^x) \to0$ outside this set].
In the case $\alpha=1$, the limit distribution
is the same as for square integrable
random walk on the plane; cf.~\cite{bf1}.
Recall that $X$ has a density $f_X$. For simplicity
we set $\gamma:=2f_X(0) {\mathbb P}(\Omega^*)$.
\begin{theorem}\label{secdistrib0}
Assume that $\alpha=1$, and fix any $x\in\RR$. Conditioned on
$\Omega^*_x$,
the variables $\varepsilon G({\tiltau_\eps^x})$
converge in law,
\[
\lim_{\eps\to0}
\PP\bigl( \gamma\eps G({\tiltau_\eps^x})
\le t | \Omega^*_x \bigr)
= \frac t {1+t} \qquad\forall t>0.
\]
\end{theorem}

For $\alpha\in(1,2]$, different limit distributions
arise, and we obtain convergence in law of $\tiltau_\eps^x$
to the $\frac1 \beta$-stable subordinator at an independent
exponential time:
\begin{theorem}\label{secdistrib}
Assume that $\alpha\in(1,2]$, and fix any $x\in\RR$.
Conditioned on $\Omega^*_x$, the variables $\varepsilon G({\tiltau
_\eps^x})$
converge in law,
\[
\lim_{\eps\to0}
\PP\biggl( {\Gamma \biggl(\frac{1}{\beta}\biggr)}\frac
\gamma\beta \eps G({\tiltau_\eps^x})
\le t \Big| \Omega^*_x \biggr)
= \Pr( {\mathcal E} {\mathcal G}_{1/\beta}^{1/\beta}
\le t ) \qquad\forall t>0
\]
or, equivalently,
\[
\lim_{\eps\to0}
\PP\biggl( \biggl({\Gamma \biggl(\frac{1}{\beta}\biggr)} \frac
\gamma\beta\biggr)^\beta
\frac{\tiltau_{\eps}^x}{G^{-1}(1/\eps)}
\le t \Big| \Omega^*_x \biggr)
= \Pr( {\mathcal E}^\beta
{\mathcal G}_{1/\beta} \le t )
\qquad\forall t>0,
\]
where ${\mathcal E}$ and ${\mathcal G}_{1/\beta}$ are
independent random variables, $\Pr({\mathcal E}>t)=e^{-t}$
and ${\mathcal G}_{1/\beta}$ having the
one-sided stable law of index $\frac{1}{\beta}$ with Laplace
transform ${\mathbb E}[e^{-s{\mathcal G}_{1/\beta}}]
= e^{-s^{1/\beta}}$, $s>0$.
\end{theorem}

In particular, we have:
\begin{coro}\label{secCOROL}
If $(X_n)_{n\ge1}$ is an i.i.d. sequence of centered
random variables with variance $1$,
satisfying the Cram\'er condition, and $x\in\RR$,
then
\[
\lim_{\eps\to0}
\PP\bigl( 2 \PP(\Omega^*) \eps\sqrt{\tiltau_\eps^x}
\le t | \Omega^*_x \bigr)
= \Pr\biggl( \frac{\mathcal E}{|\mathcal N|}
\le t \biggr)\qquad\forall t>0
\]
or, equivalently,
\[
\lim_{\eps\to0}
\PP\bigl( 4 \PP(\Omega^*)^2 \eps^2 {\tiltau_\eps^x}
\le t | \Omega^*_x \bigr)
= \Pr\biggl( \biggl(\frac{\mathcal E}{|\mathcal N|}\biggr)^2
\le t \biggr)\qquad
\forall t>0,
\]
where $\mathcal E$ and
$\mathcal N$ are independent variables,
$\mathcal N$ having a standard Gaussian distribution
${\mathcal N}(0,1)$.
\end{coro}

As Cheliotis does in~\cite{Cheliotis}, we will use the following
extension of Stone's local limit theorem~\cite{Stone}.
\begin{prop}\label{lem6}
Let $\theta$ be such that $\limsup_{\vert t\vert\rightarrow\infty}
|{\mathbb E}[e^{itX_1}]|<\theta<1$, and let $c>1$.
Then there exists a real number $h_0>0$ and an integer
$n_0\ge1$ such that, for any $n\ge n_0$, for any interval
$I$ contained in $[-h_0,h_0]$, of length larger than $\theta^n$,
we have
\[
c^{-1} f_X(0)\vert I\vert<{{\mathbb P}\biggl(\frac{S_n}{A_n}\in
I\biggr)}
<c f_X(0)\vert I\vert.
\]
\end{prop}

%
%s2 ###
\section{\texorpdfstring{Almost sure convergence: Proof of Theorem \protect\ref{thmRWP2}}{Almost sure convergence: Proof of Theorem 1}}

\mbox{}

\begin{pf*}{Proof of Theorem~\ref{thmRWP2}}
To begin with, choose $\theta$, $c$ and $h_0$ as in
Proposition~\ref{lem6}.\vadjust{\goodbreak}

To first establish an estimate from below, we
fix any $ \xi>1 $ and set $\varepsilon_n:=G(n)^{-\xi}$.
This makes the series
$\sum_n \PP(|S_n|<\varepsilon_n)$ summable:
Indeed, by regular variation and~(\ref{eqGasymptotics}),
we have $\frac{\varepsilon_n}{A_n}>\theta^n$ for $n$
large, while
\[
\frac{\eps_n}{A_n}
=O\biggl( \frac{G(n)-G(n-1)}{G(n-1)^\xi}\biggr)
=O\biggl(\int_{n-1}^n \frac{G'(t)}{G(t)^\xi} \,dt\biggr),
\]
which is summable since
$ \int_1^\infty\frac{G'(t)}{G(t)^\xi} \,dt=
[\frac{G(t)^{1-\xi}}{1-\xi}]_1^\infty<\infty$.
In particular,\break $(\frac{-\eps_n}{A_n}$, $\frac{\eps_n}{A_n})
\subseteq[-h_0,h_0]$ for large $n$. Proposition~\ref{lem6}
therefore applies to these intervals and shows that
$\PP(\vert S_n\vert<\eps_n)=O(\frac{\eps_n}{A_n})$
is summable\vspace*{1pt} as well. Hence, by the Borel--Cantelli lemma,
$\PP(\vert S_n\vert<\eps_n\mbox{ i.o.})=0$. Since
$\varepsilon_n\searrow0$, we can \mbox{conclude} that
$\tiltau_{\varepsilon_n}>n$ eventually,
almost surely on $\Omega^*$, and we get\break
$ \liminf_{n\to\infty}
\frac{\log G(\tiltau_{\varepsilon_n})}{-{\log\varepsilon_n}}\ge
\frac1\xi$
a.s. on $\Omega^*$.
Using monotonicity of $\log G(\tiltau_{\varepsilon})$
and the fact that $\eps_{n+1}\sim\eps_n$, this extends
from the $\eps_n$ to the full limit as $\eps\to0$,
and since $\xi>1$ was arbitrary, we conclude that
%
%e4 ###
\begin{equation}\label{eqliminf}
\liminf_{\epsilon\to0} \frac{\log G(\tiltau_\eps)}{-{\log
\varepsilon}} \ge1\qquad
\mbox{a.s. on } \Omega^*.
\end{equation}

To control the corresponding $\limsup$, we now fix any $\xi\in(0,1)$.
{}From Proposition~\ref{lem6}, using intervals
$(\frac{-\eps_n}{A_n},\frac{\eps_n}{A_n})$ and regular variation of
$(A_n)_{n\ge1}$, we see that there exists a constant
$c'>0$ such that for every $\varepsilon\in(0,1)$ there
is some $m_\varepsilon$ satisfying
\[
\PP(|S_k|<\varepsilon) \ge\frac{c' \varepsilon}{A_k}\qquad
\mbox{for } k \ge m_\varepsilon.
\]
More precisely, the dependence of $m_\varepsilon$ on
$\eps$ comes from the requirement $2\eps/A_k>\theta^k$ for
$k\ge m_\varepsilon$ on the length of
intervals, which is met by taking
$m_\varepsilon:=\kappa(-{\log\eps})$ with a suitable
constant $\kappa>0$.
Next, choose integers $n_\eps$ in such a way that
$G(n_\eps) \le\eps^{-1/\xi} < G(n_\eps+1)$.
Inspired by a decomposition used by
Dvoretski and Erd\"os~\cite{DvoretskiErdos},
we consider the pairwise disjoint events
$E_k^{\eps}:=\{|S_k|<\varepsilon$ and $\forall j=k+1,
\ldots,n_\eps\dvtx |S_j-S_k|>2\varepsilon\}$, $1\le k\le{n_\eps}$.
By independence and stationarity we have
\[
1 \ge\sum_{k=m_\eps}^{n_\eps} \PP(E_k^{\eps})
\ge\sum_{k=m_\eps}^{n_\eps} \PP(|S_k|<\varepsilon)
\PP(\tiltau_{2\varepsilon}>{n_\eps}-k)
\ge c'\eps \PP(\tiltau_{2\varepsilon}> {n_\eps})
\sum_{k=m_\eps}^{n_\eps} \frac{1}{A_k}.
\]
Combining this with $G(m_\eps)=o(G(n_\eps))$ [note that
$G(m_\eps)$ is slowly varying], we obtain
\begin{eqnarray*}
\PP\bigl(G(\tiltau_{2\varepsilon}) > \eps^{-{1}/{\xi}}\bigr)
&\le&\PP\bigl(G(\tiltau_{2\varepsilon}) > G(n_\eps)\bigr)
= \PP(\tiltau_{2\varepsilon}>{n_\eps})
\\
&\le&\frac{1}{c'\eps (G(n_\eps)-G(m_\eps))}
\sim\frac{\varepsilon^{1/\xi-1}}{c'}.
\end{eqnarray*}
Therefore, if we let $\varepsilon_p:=p^{-2/({1-\xi})}$,
$p\ge1$, the Borel--Cantelli lemma implies
$G(\tiltau_{2\varepsilon_p})\le
\eps_p^{-1/\xi}$ eventually almost surely,
showing that
\[
\limsup_{p\rightarrow+\infty}
\frac{\log G(\tiltau_{2\varepsilon_p})}{-{\log}(2\eps_p)}\le\frac
1\xi.
\]
Using monotonicity as before, we can extend this from
the $\eps_p$ to the full limit $\eps\to0$, and
since this is true for any $\xi\in(0,1)$, we obtain
%
%e5 ###
\begin{equation}\label{eqlimsup}
\limsup_{\eps\rightarrow0}
\frac{\log G(\tiltau_{\varepsilon})}{-{\log}(\eps)}
\le1 \qquad\mbox{a.s. on } \Omega.
\end{equation}

%%%%%%%%%%%%%%%%%%%%%%%%%%%%%%%%%%%%%%%%%%%%%%%%%%%%%%%%%%%%%%%%%%%%%%%%%%%%%%%%%%%%%%%%%%%%%%%%%%%%%%%%%%%%%%%%%%%%%%%%%
To conclude the proof, we note that for any
$\alpha\in[1,2]$ we have
\[
\lim_{n\to\infty} \frac{\log G(n)}{\log n}=\frac1\beta,
\]
which follows readily from regular variation of $G$;
compare Fact 2 in~\cite{Cheliotis}.
Together with (\ref{eqliminf}) and (\ref{eqlimsup}),
this entails
\[
\lim_{\eps\to0} \frac{\log\tiltau_{\varepsilon}}{-{\log\eps}}
=\lim_{\eps\to0} \frac{\log\tiltau_{\varepsilon}}{\log
G(\tiltau_{\varepsilon})}
\cdot
\frac{\log G(\tiltau_{\varepsilon})}{-{\log\eps}}
=\beta \qquad\mbox{a.s. on } \Omega^*
\]
as required.
\end{pf*}

The first argument can easily be adapted to prove the lower
bound (\ref{eqliminf}) also for $\tiltau_\eps^x$ with $x\neq0$.
%However, we were not able to obtain the upper bound \eqref{eqlimsup},
%even using the idea of Section~\ref{secfinal}.

%

%%%%%%%%%%%%%%%%%%%%%%%%%%%%%%%%%%%%%%%%%%%%%%%%%%%%%%%%%%%%%%%%%%%%%%%%%%%%%%%%%%%%%%%%%%%%%%%%%%%%%%%%%%%%%%%%%%%%%%%%%%%

%s3 ###
\section{Convergence in distribution for auxiliary processes}

We need to introduce auxiliary processes.
Let $(M_0^\varepsilon)_{\varepsilon>0}$
be a family of random variables, independent of
$(S_n)_{n\geq0}$, such that
$M_0^\varepsilon$ has uniform distribution on the
interval $(-\varepsilon,\varepsilon)$.
For each $\varepsilon>0$ we define
the walk $(M_n^\varepsilon)_{n\geq0}$ with random initial
position $M_0^\varepsilon$, that is,
$M_n^\varepsilon:=M_0^\varepsilon+S_n$.

A major step toward Theorems~\ref{secdistrib0}
and~\ref{secdistrib} will be to prove a version
which applies to the variables
\[
\tau_\varepsilon:=\min\{n\ge1\colon\vert
M_n^{\varepsilon}\vert
<\varepsilon\}, \qquad\eps>0.
\]
That is, we are interested in the limiting behavior, as
$\varepsilon\to0$, of the first return time distribution of the walk
$(M_n^\varepsilon)_{n\ge0}$ to the interval
$(-\varepsilon,\varepsilon)$.
The goal of the present section is to establish:
\begin{theorem}\label{thmtaueps0}
Assume that $\alpha=1$. Conditioned on $\Omega^*$,
the variables $\varepsilon G({\tau_\varepsilon})$
converge in law,
%
%e6 ###
\begin{equation}\label{eqtaueps0}
\lim_{\eps\to0}
\PP\bigl( \gamma\eps G({\tau_\varepsilon})
\le t | \Omega^* \bigr)
= \frac t {1+t} \qquad\forall t>0.
\end{equation}
\end{theorem}
\begin{theorem}\label{thmtaueps}
Assume that $\alpha\in(1,2]$.
Conditioned on $\Omega^*$, the variables $\varepsilon G({\tau
_\varepsilon})$
converge in law,
%
%e7 ###
\begin{equation}\label{eqtaueps}
\lim_{\eps\to0}
\PP\biggl( {\Gamma \biggl(\frac{1}{\beta}\biggr)}\frac
\gamma\beta \eps G({\tau_\varepsilon})
\le t \Big| \Omega^* \biggr)
= \Pr( {\mathcal E} {\mathcal G}_{1/\beta}^{1/\beta}
\le t )\qquad\forall t>0.
\end{equation}
Equivalently,
\[
\lim_{\eps\to0}
\PP\biggl( \biggl({\Gamma \biggl(\frac{1}{\beta}\biggr)} \frac
\gamma\beta\biggr)^\beta
\frac{\tau_\varepsilon}{G^{-1}(1/\eps)}
\le t \Big| \Omega^* \biggr)
= \Pr( {\mathcal E}^\beta
{\mathcal G}_{1/\beta} \le t )
\qquad\forall t>0.
\]
\end{theorem}

Again we start with considerations valid for any $\alpha\in[1,2]$.
To begin with, we define, for $\eps>0$, $R>0$, and integers $K>0$,
auxiliary events
\[
\Gamma_{\eps,R,K}:=
\{\forall i=1,\ldots,K\dvtx S_i\neq0 \mbox{ and } |M^\eps_i|\le R\},
\]
which asymptotically exhaust $\Omega^*$, and on which we can
work conveniently.
As $\eps\to0$ we have
$\PP(\Gamma_{\eps,R,K}) \to\PP(\Gamma_{R,K})$ and
$\PP(\Gamma_{\eps,R,K}\setminus\Omega^*)
\to\PP(\Gamma_{R,K} \setminus\Omega^*)$,
where $\Gamma_{R,K} := \{\forall i=1,\ldots,K \dvtx 0<|S_i|\le R\}$
(except, perhaps, for a countable set of $R$'s which we are
going to avoid). Let $n\in\mathbb{N}$.
Using again a decomposition similar to that
of Dvoretski and Erd\"os in
\cite{DvoretskiErdos}, we find, for $\eps\in(0,\frac1 2)$,
%
%e8 ###
\begin{equation}\label{eqdvoer}
\PP(\Gamma_{\eps,R,K}) = \sum_{k=0}^n p_{k}^- = \sum_{k=0}^n p_{k}^+
\end{equation}
with
$p_{k}^\pm=p_{k,n,\eps,R,K}^\pm:= \PP(\Gamma_{\eps,R,K} \cap\{
|M^\eps_k|<\eps\pm2 \eps^2
\mbox{ and } \forall\ell=k+1,\ldots,n\dvtx \break|M^\eps_\ell|\ge\eps\pm
2 \eps^2\})$ for
$1 \le k \le n$,
and $p_{0}^\pm=p_{0,n,\eps,R,K}^\pm:=\PP(\Gamma_{\eps,R,K} \cap\{
\forall\ell=1,\ldots,n \dvtx |M^\eps_\ell|\ge\eps\pm2 \eps^2\})$.
In the sequel, we will use the following notation:
given two functions $a$ and $b$, the notation
\[
a(\varepsilon,R,K)=o_{\varepsilon,R,K}(1)
\quad\mbox{and}\quad b(R,K)=o_{R,K}(1)
\]
will mean that
\[
\limsup_{K\rightarrow+\infty}\limsup_{R\rightarrow+\infty}\limsup
_{\varepsilon\rightarrow0}
|a(\varepsilon,R,K)|=0 \quad\mbox{and}\quad
\limsup_{K\rightarrow+\infty}\limsup_{R\rightarrow+\infty} |b(R,K)|=0.
\]
We will also write $m_\varepsilon:=(\log\varepsilon)^4$.
The following estimates are the basis of the argument to follow.
\begin{lem}\label{seclemupper}
Let $c>0$ and let, for every $\varepsilon>0$, $n_\varepsilon$ be the
integer such that
$G(n_\varepsilon)\le\frac c\varepsilon<G(n_\varepsilon+1)$.

For arbitrary $\gamma'$ and $\gamma''$ such that $0<\gamma
'<2f_X(0)<\gamma''$, we have
\[
\PP(\Gamma_{\eps,R,K})
\ge
\PP(\tau_\eps>n_\eps)
+{\PP(\Gamma_{\eps,R,K}) {\gamma'}\eps}
\sum_{k=m_\eps}^{n_\eps}
\frac{\PP(\tau_\varepsilon>n_\eps-k)}{A_{k}}+o_{\varepsilon,R,K}(1)
\]
and
\[
\PP(\Gamma_{\eps,R,K})
\le
\PP(\tau_\eps>n_\eps) + \PP(\Gamma_{\eps,R,K}) {\gamma''}
\eps\sum_{k=m_\eps}^{n_\eps} \frac{\PP(\tau_\eps>n_\eps-k)}{A_k}
+o_{\varepsilon,R,K}(1).
\]
\end{lem}
\begin{pf}
For the course of this proof, we simplify notation by suppressing
the parameters $\eps$, $R$, and $K$ in $m_\eps$, $n_\eps$, $M^\eps
_i$, and $\Gamma_{\eps,R,K}$.
We will apply (\ref{eqdvoer}) with $n=n_\eps$.
Also, let $\nu:=\eps^2$.

\begin{longlist}
\item
Starting with the $k=0$ term, we see that
\[
p_{0}^{-}\ge\PP(\Gamma\cap\{\forall\ell=1,\ldots,n:
|M_\ell|\ge\eps\})
\ge
\PP(\Gamma\cap\{\tau_\eps>n\}).
\]
We now consider the case where $m\le k \le n$.
Let $\A:=(2\nu\ZZ)\cap(-\eps+3\nu,\eps-3\nu)$.
Notice that the sets $Q_a:=(a-\nu,a+
\nu)$ with $a\in\A$ are disjoint and
contained in $(-\eps+2\nu,\eps-2\nu)$.
Therefore the $k$th term in (\ref{eqdvoer}) satisfies
%
%e9 ###
\begin{eqnarray}\label{eqN1}\quad
p_{k}^{-}&\ge&
\sum_{a\in\A}\PP(\Gamma\cap\{
M_k\in Q_a \mbox{ and } \forall\ell=k+1,\ldots,n\dvtx
\vert M_\ell\vert\ge\eps-2\nu\})\nonumber\\
&\ge&
\sum_{a\in\A}\PP(\Gamma\cap\{
M_k\in Q_a \mbox{ and } \forall\ell=k+1,\ldots,n\dvtx
\vert S_\ell-S_k+a\vert\ge\eps-\nu\}) \\
&=&
\sum_{a\in\A}\PP(\Gamma\cap\{
M_k\in Q_a\})
\PP( \forall\ell=1,\ldots,n-k\dvtx
\vert S_\ell+a\vert\ge\eps-\nu)\nonumber
\end{eqnarray}
by independence [where we assume that $\eps$ is so small that
$(\log\eps)^4>K$]. Note that
\begin{eqnarray*}
&&\PP(\Gamma\cap\{ M_k\in Q_a\}) \\
&&\qquad= \int_{\{
\forall i\dvtx x_i\neq x_0,|x_i|\le R\}} \PP(S_{k-K}\in Q_a-x_K)\, d\PP
_{(M_0,\ldots,M_K)}(x_0,\ldots,x_K)
\end{eqnarray*}
with $d\PP_{(M_0,\ldots,M_K)}$ denoting the distribution of
$(M_0,\ldots,M_K)$.
Now fix $\theta$ as in Proposition~\ref{lem6}, and $c\in(0,1)$ such that
$\gamma' <2f_X(0)/c$. Elementary considerations
show that Proposition~\ref{lem6} applies to $I=\frac1{A_{k-K}}(Q_a-x_K)$
if $\eps$ is sufficiently small, and in this case gives
%
%e10 ###
\begin{equation}\label{eqN2}
\PP(\Gamma\cap\{M_k\in Q_a\}) \ge\PP(\Gamma)
\frac{\gamma'\nu}{A_{k}}.
\end{equation}
Using this, plus the observation that conditioning on $\{M_0\in Q_a\}$
amounts to looking at $M^*_n:=M^*_0+S_n$, $n\ge0$, with $M^*_0$
uniformly distributed on $Q_a$,
we can continue to estimate, for small $ \eps$,
%
%e11 ###
\begin{eqnarray}\label{eqN3}\qquad
p_{k}^{-}&\ge&
\PP(\Gamma) \frac{\gamma'\nu}{A_{k}} \sum_{a\in\A}
\PP( \forall\ell=1,\ldots,n-k\dvtx \vert S_\ell+a\vert\ge\eps-\nu
) \nonumber\\
&\ge&
\PP(\Gamma) \frac{\gamma'\nu}{A_{k}} \sum_{a\in\A}
\PP( \{\forall\ell=1,
\ldots,n-k\dvtx |M_\ell| \ge\eps\} | \{M_0\in Q_a\})\nonumber\\
&\ge&
\PP(\Gamma) \frac{\gamma'\eps}{A_{k}} \sum_{a\in\A}
\PP( \{\forall\ell=1,
\ldots,n-k\dvtx |M_\ell| \ge\eps\} \cap\{M_0\in Q_a\})\\
&\ge&
\PP(\Gamma)
\frac{\gamma'\eps}{A_{k}}
\bigl( \PP( \forall\ell=1,\ldots,n-k\dvtx
|M_\ell| \ge\eps)-
\PP(\eps-4\nu\le|M_0|\le\eps) \bigr)\nonumber\\
&=&
\PP(\Gamma)
\frac{\gamma'\eps}{A_{k}}
\bigl( \PP( \tau_\eps>n-k)- 8\nu\bigr).\nonumber
\end{eqnarray}
Putting together these estimates via equation (\ref{eqdvoer}) gives
\begin{eqnarray*}
&&\PP(\Gamma\cap\{\tau_\eps>n\})
+{\PP(\Gamma) \gamma' \eps}
\sum_{k=m}^{n}
\frac{\PP(\tau_\varepsilon>n-k)}{A_{k}}\\
&&\qquad\le
\PP(\Gamma) + \PP(\Gamma)8\gamma'\eps\nu\bigl(G(n)-G(m)\bigr).
\end{eqnarray*}
We observe that $\Gamma^c \cap\{\tau_\eps>n\} \subseteq
\bigcup_{i=1}^K\{\vert M_i\vert>R\}$ for $\eps$
so small that $n=n_\eps> K$.
Since ${\limsup}_{K\rightarrow+\infty}{\limsup}_{R\rightarrow
+\infty}{\limsup}_{\eps\rightarrow0}
{\mathbb P}(\bigcup_{i=1}^K\{\vert M_i\vert>R\})=0$
and $\lim_{\eps\rightarrow0}\eps^3(G(n)-G(m))=0$,
this proves the first assertion of the lemma.

\item
We only provide a sketch of the proof
of the second point since the arguments
are very similar to the above.
Using (\ref{eqdvoer}) gives
\[
\PP(\Gamma)\le
\PP(\Gamma\cap\{\tau_\eps>n\})+
\PP(\Gamma\setminus\Omega^*) + \PP(\Omega^* \cap\{\tau_{3\eps
}\le m\})+\sum_{k=m}^{n} p_k^+
\]
since $ \sum_{k=1}^{m} p_k^+ \le\PP(\Gamma\cap
\{\tau_{3\eps}\le m\})$.
Next, take $\bar\A:=(2\nu\ZZ)
\cap(-\eps-3\nu,\eps+3\nu)$ and intervals $\bar Q_a:=[a-\nu,a+
\nu]$, $a\in\bar\A$, which cover $(-\eps-2\nu,\eps+2\nu)$.
We can then use arguments parallel to those of part (i)
to obtain
\begin{eqnarray*}
\sum_{k=m}^{n}p_k^+ &\le&
\sum_{k=m}^{n}\sum_{a\in\bar\A}
\PP\bigl(\Gamma\cap\{M_k\in\bar Q_a \mbox{ and }
\forall\ell=k+1,\ldots,n \dvtx \vert M_\ell\vert>
\varepsilon+2\nu\}\bigr)\\
& \vdots&\\
&\le&
\PP(\Gamma) {\gamma''} \eps\sum_{k=m}^{n}
\frac{\PP(\tau_\eps>n-k)}{A_k}
+ \PP(\Gamma)8\gamma''\eps\nu\bigl(G(n)-G(m)\bigr),
\end{eqnarray*}
which\vspace*{1pt} proves our claim since
$\lim_{\eps\rightarrow0}\PP(\Omega^* \cap\{\tau_{3\eps}\le m\}
)=0$ as a consequence of Theorem
\ref{thmRWP2} and since
$\PP(\Gamma\setminus\Omega^*)=o_{\eps,R,K}(1)$.\qed
\end{longlist}
\noqed\end{pf}

This enables us to
derive an asymptotic bound for the tails of the
distributions of the $\eps G(\tau_\eps)$ as $\eps\to0$.
\begin{lem}\label{seclimsup}
For all $\alpha\in[1,2]$ and any $t>0$ we have
\[
\limsup_{\eps\to0}
\PP\bigl(\gamma\eps
G(\tau_\eps)>t\bigr)\le
\frac{\PP(\Omega^*)}
{1+{t} }.\vspace*{-2pt}
\]
\end{lem}
\begin{pf}
Fix $t$, $R$, $K$ and $0<\gamma'<2f_X(0)$.
For $\eps>0$ choose $n_\eps$ so that
$G(n_\eps)\le\frac t{\gamma\eps}\le G(n_\eps+1)$, whence
$\PP(\eps\gamma G(\tau_\eps)>t) \sim\PP(\tau_\eps>n_\eps)$.
Recall that $m_\eps:=(\log\eps)^4$.
As in the proof of Theorem~\ref{thmRWP2} we see that
$G(m_\eps)=o(G(n_\eps))$. Therefore
%
%e12 ###
\begin{equation}\label{eqM1}\qquad
\eps
\sum_{k=m_\eps}^{n_\eps}
\frac{\PP(\tau_\eps>n_\eps-k)}{A_{k}}
\ge\eps\bigl(G(n_\eps)-G(m_\eps)\bigr)\PP(\tau_\eps>n_\eps)
\sim
\frac{t}{\gamma}
\PP(\tau_\eps>n_\eps).
\end{equation}
Together with the first part of Lemma~\ref{seclemupper},
this yields
\[
\limsup_{\eps\to0} \PP\bigl(\eps\gamma
G(\tau_\eps)>t\bigr)\le
\frac{\PP(\Gamma_{R,K})+ o_{R,K}(1)}
{1+(t{\gamma'}/{\gamma}) \PP(\Gamma_{R,K}) }.
\]
Taking successively $R\to\infty$, then
$K\to\infty$ and finally $\gamma'\to2f_X(0)$,
we obtain the lemma.\vspace*{-2pt}
\end{pf}

When $\alpha=1$, this upper bound actually is the limit:\vspace*{-2pt}
\begin{lem}\label{secalpha=1}
If $\alpha=1$, then for any $t>0$ we have
\[
\liminf_{\eps\to0} \PP\bigl(\gamma\eps
G(\tau_\eps)>t\bigr)\ge
\frac{\PP(\Omega^*)}
{1+{t}}.\vspace*{-2pt}
\]
\end{lem}
\begin{pf}
Fix $t$, $R$, $K$ and $\gamma''>2f_X(0)$,
and choose $m_\eps$ and $n_\eps$
as in the previous proof.

Since $\alpha=1$ means that $G$ is slowly varying,
we have $G(2n_\eps)-G(n_\eps)=o(G(n_\eps))$. Hence
%
%e13 ###
\begin{eqnarray}\label{eqM2}
%&\le
&&\PP(\tau_\eps>2n_\eps) +
{\PP(\Gamma_{\eps,R,K}) {\gamma''}\eps}
\sum_{k=m_\eps}^{2n_\eps}
\frac{\PP(\tau_\eps>2n_\eps-k)}{A_{k}}\nonumber\\[-2pt]
&&\qquad\le
\PP(\tau_\eps>n_\eps)+
{\PP(\Gamma_{\eps,R,K}) {\gamma''}\eps}
\Biggl(
\sum_{k=m_\eps}^{n_\eps}
\frac{\PP(\tau_\eps>n_\eps)}{A_{k}}
+\sum_{k=n_\eps}^{2n_\eps}
\frac{1}{A_{k}}
\Biggr) \nonumber\\[-9pt]\\[-9pt]
&&\qquad\le
\PP(\tau_\eps>n_\eps)+
{\PP(\Gamma_{\eps,R,K}) {\gamma''}\eps}G(n_\eps)
[{\PP(\tau_\eps>n_\eps)}+o(1)]\nonumber\\[-2pt]
&&\qquad\le
\PP(\tau_\eps>n_\eps)+
t\frac{\gamma''}{\gamma} \PP(\Gamma_{\eps,R,K})
{\PP(\tau_\eps>n_\eps)}+o(1).\nonumber
\end{eqnarray}
Combining these observations with
the second estimate of Lemma~\ref{seclemupper}
(replacing $n_\eps$ by $2n_\eps$) entails
\[
\liminf_{\eps\rightarrow0}\PP(\tau_\eps>n_\eps)
\ge
\frac{\PP(\Gamma_{R,K})-o_{R,K}(1)}
{1+(t{\gamma''}/{\gamma}) \PP(\Gamma_{R,K})}.
\]
We conclude by successively taking
$R\rightarrow\infty$, $K\rightarrow\infty$ and
$\gamma''\rightarrow2f_X(0)$.\vadjust{\goodbreak}
\end{pf}

\begin{pf*}{Proof of Theorem~\ref{thmtaueps0}}
Immediate from Lemmas~\ref{seclimsup} and~\ref{secalpha=1},
as $\eps G(\tau_\eps) \to0$ outside $\Omega^*$.
\end{pf*}

When $\alpha\in(1,2]$,
Lemma~\ref{seclemupper}
does not yet give the limit distribution. Still, it immediately implies
the tightness of the family
of distributions with the normalization given there:
\begin{lem}\label{lemtight}
The family of distributions of the random variables $\eps G(\tau_\eps
)$, $\eps\in(0,1)$, is tight.
\end{lem}

Hence it will be enough to prove that the advertised limit law
is the only possible accumulation point of our distributions.
We henceforth abbreviate
\[
Z_\eps:= \frac\gamma\beta \eps G({\tau_{\eps}}),\qquad \eps>0.
\]

\begin{lem}\label{seclower}
Suppose that $\alpha\in(1,2]$.
Let $(\eps_p)_{p\ge1}$ be a positive sequence with
$\lim_{p\to\infty}\eps_p=0$, and such that the conditional
distributions of the $Z_{\eps_p}$ on $\Omega^*$ converge to the law
of some
random variable $Y$. Then its tail satisfies the integral equation
\[
1 = {\Pr}(Y>t)+t\int_0^1
\frac{\Pr(Y>t(1-u)^{1/\beta})}{u^{1/\alpha}}\,
du \qquad\forall t>0.
\]
\end{lem}
\begin{pf}
(i)
We write $f(t):= {\Pr}(Y>t)$, and first prove that
\[
\forall t>0\qquad
1 \ge f(t)+t\int_0^1u^{-1/\alpha}
f\bigl(t(1-u)^{1/\beta}\bigr) \,du.
\]
Let us only consider $\varepsilon$ belonging to
$\{\eps_p, p\ge1\}$.
Note that by monotonicity and right continuity of $f$ it suffices
to prove the inequality for all
$t\in(0,\infty)$ such that, for all
$N\ge1$ and all $r=0,\ldots,N-1$,
the function $f$ is continuous at
$t(1-\frac{r}{N})^{1/\beta}$.
Henceforth such a $t$ will be fixed.

Now take some $\delta>0$. We claim that one can choose $N_\delta> 1$
such that for all $N \ge N_\delta$,
%
%e14 ###
\begin{equation}\label{eqRI}
\Biggl\vert\int_0^1\frac{f(t(1-u)^{1/\beta})}
{{u}^{1/\alpha}} \,du
-\frac{1}{N}\sum_{r=1}^{N-1}
\frac{f(t(1-(r/N))^{1/\beta})}
{((r+1)/N)^{1/\alpha}}
\Biggr\vert\le\delta.
\end{equation}
Indeed, take $\Delta\in(0,1)$ such that $\beta\Delta^{1/\beta}
<\delta/4$.
For any $N$ we have
\[
\frac1N \sum_{r=1}^{\lfloor\Delta N\rfloor} \biggl(\frac
{r+1}{N}\biggr)^{-1/\alpha}
\le\int_0^\Delta u^{-1/\alpha} \,du= \beta\Delta^{1/
\beta} < \delta/4.
\]
Since $f$ is bounbed by one this implies that both the integral in
(\ref{eqRI})
restricted to $[0,\Delta]$ and the sum from $r=1$ to $\lfloor\Delta
N\rfloor$
are bounded by $\delta/4$.
The claim follows by taking $N_\delta$ so large that the approximation
of the Riemann
integral on the interval $[\Delta,1]$ by the Riemann sum with step
$1/N$ has a precision
at least $\delta/2$.

Now fix integers $N \ge N_\delta$,
$K\ge1$, and some $0<\gamma'<2f_X(0)$.
For $\eps>0$ small enough take $n_\eps$ such that
$G(n_\eps)\le\frac{\beta t}{\gamma\eps}<G(n_\eps+1)$
[and hence \mbox{$G(n_\eps)\sim\frac{\beta t}{\gamma\eps}$}].

According to the first point of
Lemma~\ref{seclemupper}, since $\frac{n_\eps}N\ge m_\eps$, we have
\[
\PP(\Gamma_{\eps,R,K})
\ge
\PP(Z_\eps>t)
+{\PP(\Gamma_{\eps,R,K}) {\gamma'}\eps}
\sum_{k=n_\eps/N}^{n_\eps}
\frac{\PP(\tau_\varepsilon>n_\eps-k)}{A_{k}}
+o_{\eps,R,K}(1).
\]
Due to our assumption on the $Z_{\eps_p}$ and $t$, we see that
$\PP(Z_\eps>t) \to\PP(\Omega^*) f(t)$ as $\eps_p\to0$.
Next, by monotonicity,
\begin{eqnarray*}
&&\sum_{k=n_\eps/N}^{n_\eps}
\frac{\PP(\tau_\varepsilon>{n_\eps}-k) } {A_{k}}\\
&&\qquad\ge
\sum_{r=1}^{N-1} \sum_{k=0}^{{n_\eps}/N-1}
\frac{\PP(\tau_\varepsilon>{n_\eps}-k-(r {n_\eps}/N))}
{A_{k+(r {n_\eps}/N)}}\\
&&\qquad\ge
\sum_{r=1}^{N-1} \biggl( G\biggl(\frac{r+1}{N} {n_\eps}
\biggr)-G\biggl(\frac{r}{N} {n_\eps}\biggr)\biggr)
\PP\biggl(\tau_\varepsilon>\biggl(1-\frac{r}{N}\biggr){n_\eps
}\biggr).
\end{eqnarray*}
By regular variation, the first term of the product is asymptotically
equivalent to
\[
G({n_\eps})\biggl[\biggl(\frac{r+1}{N}\biggr)^{1/\beta}-
\biggl(\frac{r}{N}\biggr)
^{1/\beta}\biggr]
\ge
\frac{G({n_\eps})}{\beta N (({r+1})/{N})^{
{1}/{\alpha}}}
\]
as $\eps_p\to0$. On the other hand, the second term is equal to
\[
\PP\biggl(Z_\eps>\eps\frac\gamma\beta G\biggl(\biggl(1-\frac r
N\biggr){n_\eps}\biggr)\biggr)
\to
\PP(\Omega^*) f\biggl( t \biggl(1-\frac r N\biggr)^{1/\beta}
\biggr),
\]
since $G((1- \frac r N){n_\eps})
\sim(1- \frac r N)
^{1/\beta}G({n_\eps})$.
As a consequence, we see that
%
%e15 ###
\begin{eqnarray}\label{eqO}
&&
\liminf_{p\to\infty}
\eps_p\sum_{k=n_{\eps_p}/N}^{n_{\eps_p}}
\frac{\PP(\tau_{\varepsilon_p}>{n_{\eps_p}}-k) } {A_{k}}\nonumber\\
&&\qquad\ge
\PP(\Omega^*) \frac{t}{\gamma} \frac1 N
\sum_{r=1}^{N-1}
\frac{
f(t (1- r/ N)^{1-1/\alpha})}
{(({r+1})/N)^{1/\alpha}}\\
&&\qquad\ge
\PP(\Omega^*) \frac{t}{\gamma} \biggl(
\int_0^1 \frac{f(t(1-u)^{1-1/\alpha})} {{u}^{1/\alpha}}\,
du -\delta
\biggr).
\nonumber
\end{eqnarray}
Combining all these asymptotic estimates and taking the limit
$\eps_p\to0$, we end then up with
\begin{eqnarray*}
\PP(\Gamma_{R,K})
&\ge&
\PP(\Omega^*) \biggl[ f(t)
+\frac{\PP(\Gamma_{R,K}) \gamma't}{ \gamma}
\biggl(\int_0^1
\frac{f(t(1-u)^{1-1/\alpha})} {{u}^{1/\alpha}} \,du
-\delta\biggr) \biggr]\\
&&{}+o_{R,K}(1).
\end{eqnarray*}
Successively letting $R\to\infty$, $K\to\infty$,
$\gamma'\to2 f_X(0)$ and $\delta\to0$ we obtain the
desired inequality.

(ii) The converse inequality is proved
analogously, using the other half of
Lemma~\ref{seclemupper} with the following adaptation:
we have
\[
\PP(\Gamma_{\eps,R,K})
\le
\PP(Z_\eps>t)
+{\PP(\Gamma_{\eps,R,K}) {\gamma''}\eps}
\sum_{k=m_\eps}^{n_\eps}
\frac{\PP(\tau_\varepsilon>n_\eps-k)}{A_{k}}
+o_{\eps,R,K}(1).
\]
Since, $G(n_\eps/N)\sim G(n_\eps)N^{-1/\beta}$ as $\eps$ goes
to 0, we have,
for $\eps$ small enough,
\[
\eps\sum_{k=m_\eps}^{n_\eps/N}
\frac{\PP(\tau_\varepsilon>n_\eps-k)}{A_{k}}
\le\eps G\biggl(\frac{n_\eps}N\biggr)
\le2\eps G(n_\eps)N^{-1/\beta}\le
2 \frac{\beta t}\gamma N^{-1/\beta}
\]
and so
\begin{eqnarray*}
\PP(\Gamma_{\eps,R,K})
&\le&\PP(Z_\eps>t)
+{\PP(\Gamma_{\eps,R,K}) {\gamma''}\eps}
\sum_{k=n_\eps/N}^{n_\eps}
\frac{\PP(\tau_\varepsilon>n_\eps-k)}{A_{k}}
+o_{\eps,R,K}(1)\\
&&{}+ 2 \gamma''\frac{\beta t}\gamma N^{-1/\beta}.
\end{eqnarray*}
\upqed\end{pf}

Now let us identify the limit distribution
satisfying the equality given by Lem\-ma~\ref{seclower}.
To this end we consider the variables
\[
Z_{\eps}':= \biggl(\frac\gamma\beta\biggr)^\beta
\frac{\tau_{\eps}}{G^{-1}(1/\eps)},\qquad \eps>0.
\]

\begin{lem}\label{lemfunctional2}
The conditional distributions of the $Z_{\eps_p}$ converge
to a random variable $Y$ iff the conditional distributions
of the $Z_{\eps_p}'$ converge to $Y^\beta$. The latter
then satisfies
\[
1=\Pr(Y^\beta>t)+\int_0^t
\frac{\Pr(Y^\beta>t-v)}{v^{1/\alpha}} \,dv
\qquad\forall t>0.
\]
\end{lem}
\begin{pf}
The equivalence of the two conditional distributional convergence
statements follows from regular variation of $G^{-1}$;
see, for example, Lemma 1 of~\cite{BZ}. Suppose that they hold.\vadjust{\goodbreak} Then,
according to Lemma~\ref{seclower}, for any $t>0$, we have
\[
1 = \Pr(Y^\beta>t)+
t^{1/\beta}\int_0^1\frac{
\Pr(Y^\beta>t(1-u))}
{u^{1/\alpha}}
\,du,
\]
and the conclusion follows by a change of variables, $v=tu$.
\end{pf}
\begin{lem}\label{lemcarac}
Let $W$ be a random variable with values
in $[0,\infty)$ satisfying
%
%e16 ###
\begin{equation}\label{eqint}
\Pr(W\le t)=\int_0^t\frac{\Pr(W>t-v)}
{v^{1/\alpha}} \,dv \qquad\forall t>0.
\end{equation}
Then
\[
{\mathbb E}[e^{-sW}]=\frac1{1+c_\beta s^{1/\beta}}
\qquad\forall s>0
\]
with $c_\beta:={\Gamma(\frac{1}{\beta})}^{-1}$.
In particular, the distribution of $W$ coincides with that
of $c_\beta^\beta{\mathcal E}^\beta{\mathcal G}_{{1}/{\beta}}$,
where the independent variables
${\mathcal E}$ and ${\mathcal G}_{{1}/{\beta}}$ are as in the
statement of Theorem~\ref{secdistrib}.
\end{lem}
\begin{pf}
Let $s>0$. We have
\begin{eqnarray*}
{\mathbb E}[e^{-sW}]&=&\int_0^{+\infty}
\Pr(e^{-sW}\ge u) \,du\\
&=&\int_0^{+\infty}\Pr\biggl(W\le-\frac{\log(u)}{s}
\biggr) \,du\\
&=& \int_0^{+\infty}
\Pr(W\le v)se^{-sv} \,dv.
\end{eqnarray*}
Hence, for any $s>0$, we find
\begin{eqnarray*}
{\mathbb E}[e^{-sW}]
&=&\int_0^{+\infty}\biggl[\int_0^v
\frac{\Pr(W\ge v-w)}{w^{1/\alpha}} \,dw\biggr]
se^{-sv} \,dv\\
&=&\int_0^{+\infty}\frac1{ w^{1/\alpha} }
\biggl[\int_w^{+\infty}
{\Pr(W\ge v-w)}se^{-sv} \,dv\biggr] \,dw\\
&=&\int_0^{+\infty}\frac{e^{-sw}}{ w^{1/\alpha} }
\biggl[\int_0^{+\infty}
{\Pr(W\ge z)}se^{-sz} \,dz\biggr] \,dw\\
&=&\int_0^{+\infty}\frac{e^{-sw}}{w^{1/\alpha}}
\biggl[1-\int_0^{+\infty}
{\Pr(W\le z)}se^{-sz} \,dz\biggr] \,dw\\
&=&\int_0^{+\infty}\frac{e^{-sw}}{w^{1/\alpha}}
\,dw \cdot\bigl[1- {\mathbb E}[e^{-sW}] \bigr],
\end{eqnarray*}
and our claim about the Laplace transform of $W$ follows since
\[
\int_0^{+\infty}\frac{e^{-sw}}{w^{1/\alpha}}
\,dw
=\frac\beta{ s^{1/\beta} }\int_0^{+\infty}
e^{-z^\beta} \,dz=\frac1{c_\beta s^{1/\beta}}
\qquad\mbox{with }
c_\beta:=\frac{1}{\Gamma({1}/{\beta})}.\vadjust{\goodbreak}
\]
Given this, a routine calculation (cf. Problem XIII.11.10 of~\cite{Feller})
shows that $W$ indeed has the same Laplace transform as
$c_\beta^\beta{\mathcal E}^\beta{\mathcal G}_{{1}/{\beta}}$.
\end{pf}

\begin{pf*}{Proof of Theorem~\ref{thmtaueps}}
According to Lemma~\ref{lemtight} the family of distributions
of the $Z_\eps$, $\eps\in(0,1)$, is tight.
By Lemmas~\ref{seclower},~\ref{lemfunctional2} and~\ref{lemcarac},
the law of $c_\beta{\mathcal E} {\mathcal G}_{1/\beta}^{1/\beta}$
is the only possible accumulation point of these distributions.
\end{pf*}

%
%s4 ###
\section{\texorpdfstring{Convergence in distribution for $\tiltau_\eps^x$}
{Convergence in distribution for tau epsilon x}}\label{secfinal}

To complete the proof of Theorems~\ref{secdistrib0} and~\ref{secdistrib}
we now utilize Theorems~\ref{thmtaueps0}
and~\ref{thmtaueps}. Note first that it suffices to prove Theorems
\ref{secdistrib0} and~\ref{secdistrib} under the additional assumption
that $S'_0=0$, in which case
\[
\tiltau_\eps^x = \hat{\tiltau}_\eps^x
:= \inf\{n\ge1\dvtx |S_n-x|<\eps\}
\quad\mbox{and}\quad
\Omega^*_{x} = \hat{\Omega}^*_{x}
:= \{ S_n \neq x\ \forall n \}.
\]
Indeed, in the situation of Theorem~\ref{secdistrib0}, with arbitrary
distribution $P$ of $S'_0$, we then have
\[
\PP\bigl( \gamma\eps G({\tiltau_\eps^{x}})\le t \bigr)
=
\int_{{\mathbb R}}
\PP\bigl( \gamma\eps G({\hat{\tiltau}_\eps^{x-y}})\le t
\bigr)\,dP(y)\to
\int_{{\mathbb R}}
\PP( \hat{\Omega}^*_{x-y} ) \,dP(y)
\cdot\frac{t}{1+t}
\]
by the $P=\delta_0$ case of Theorem~\ref{secdistrib0} and dominated
convergence and analogously for Theorem~\ref{secdistrib}.

Therefore, for the remainder of this section we assume that $S'_0=0$.

Next, we observe that our key lemma (Lemma~\ref{seclemupper})
can be adapted as follows.
Let $\Gamma_{R,K}^x$ be the event defined by
\[
\Gamma_{R,K}^x:=
\{\forall i=1,\ldots,K\dvtx S_i\neq x \mbox{ and } |S_i|\le R\}.
\]

\begin{lem}\label{seclemupperbis}
Let $c>0$, and let, for every $\varepsilon>0$, $n_\varepsilon$ be the
integer such that
$G(n_\varepsilon)\le\frac c\varepsilon<G(n_\varepsilon+1)$.

For arbitrary $\gamma'$ and $\gamma''$ such that $0<\gamma
'<2f_X(0)<\gamma''$ we have
\[
\PP(\Gamma_{R,K}^x)
\ge
\PP(\tiltau_\eps^x>n_\eps)
+{\PP(\Gamma_{R,K}^x) {\gamma'}\eps}
\sum_{k=m_\eps}^{n_\eps}
\frac{\PP(\tau_\varepsilon>n_\eps-k)}{A_{k}}+o_{\eps,R,K}(1)
\]
and
\[
\PP(\Gamma_{R,K}^x)
\le
\PP(\tiltau_\eps^x>n_\eps) + \PP(\Gamma_{R,K}^x) {\gamma''}
\eps\sum_{k=m_\eps}^{n_\eps} \frac{\PP(\tau_\eps>n_\eps-k)}{A_k}
+ o_{\eps,R,K}(1).
\]
\end{lem}
\begin{pf}%[Proof of Lemma~\ref{seclemupperbis}]
We have the following analog of formula (\ref{eqdvoer}):
%
%e17 ###
\begin{equation}\label{eqdvoerbis}
\PP(\Gamma_{R,K}^x) = \sum_{k=0}^{n_\eps}p_k^{x,-} = \sum
_{k=0}^{n_\eps}p_k^{x,+}
\end{equation}
with
\[
p_0^{x,\pm}:= \PP(\Gamma_{R,K}^x \cap\{\forall\ell=1,\ldots
,n_\eps\dvtx
|S_\ell-x|\ge\eps\pm2 \eps^2\})
\]
and
\[
p_k^{x,\pm}:=
\PP(\Gamma_{R,K}^x \cap\{|S_k-x|<\eps\pm2 \eps^2
\mbox{ and } \forall\ell=k+1,\ldots,n_\eps\dvtx |S_\ell-x|\ge\eps
\pm2 \eps^2\}).
\]
We follow
the proof of Lemma~\ref{seclemupper}.

(i)
Observe first that
\[
p_{0}^{x,-}\ge
\PP(\Gamma_{R,K}^x \cap\{\tiltau_\eps^x>n_\varepsilon\}).
\]
Now consider indices with $m_\eps\le k \le n_\eps$.
With the same set $\mathcal A$ as in the proof of Lemma~\ref{seclemupper},
we find, arguing as in (\ref{eqN1}), that
\begin{eqnarray*}
p_{k}^{x,-}&\ge&
\sum_{a\in\A}\PP(\Gamma_{R,K}^x\cap\{
S_k-x\in Q_a \mbox{ and } \forall\ell=k+1,\ldots,n_\eps\dvtx
\vert S_\ell-x\vert\ge\eps-2\nu\})\\
&\ge&
\sum_{a\in\A}\PP(\Gamma_{R,K}^x\cap\{
S_k-x\in Q_a\})
\PP( \forall\ell=1,\ldots,n_\eps-k\dvtx
\vert S_\ell+a\vert\ge\eps-\nu).
\end{eqnarray*}
A proof parallel to that of (\ref{eqN2}) shows that
\[
\PP(\Gamma_{R,K}^x\cap\{S_k-x\in Q_a\}) \ge\PP(\Gamma_{R,K}^x)
\frac{\gamma'\nu}{A_{k}},
\]
if $\eps$ is sufficiently small. Therefore,
\begin{eqnarray*}
p_{k}^{x,-}&\ge&
\PP(\Gamma_{R,K}^x) \frac{\gamma'\nu}{A_{k}} \sum_{a\in\A}
\PP( \forall\ell=1,\ldots,n_\eps-k\dvtx \vert S_\ell+a\vert\ge\eps
-\nu\}) \\
&\ge&
\PP(\Gamma_{R,K}^x)
\frac{\gamma'\eps}{A_{k}}
\bigl( \PP( \tau_\eps>n_\eps-k)- 8\nu\bigr),
\end{eqnarray*}
where the second step uses an estimate contained in (\ref{eqN3}).
Continuing as in the proof of Lemma~\ref{seclemupper}, we
obtain the first assertion of our lemma.

(ii) Similar adaptations give the second assertion of the lemma.
\end{pf}

We can now complete the proofs of our
main distributional limit theorems:

\begin{pf*}{Proof of Theorem~\ref{secdistrib0}}
We go back to Lemmas~\ref{seclimsup} and~\ref{secalpha=1},
observing that we already have (\ref{eqtaueps0}) at our disposal.
Take $t\in(0,\infty)$, $R,K \ge1$ and $\gamma'<2f_X(0)<\gamma''$.
For $\eps>0$ let $m_\eps:=(\log\eps)^4$ and choose $n_\eps$, such that
$G(n_\eps)\le\frac t{\gamma\eps}\le G(n_\eps+1)$, meaning that
$\PP(\eps\gamma G(\tiltau_\eps^x)>t) \sim\PP(\tiltau_\eps
^x>n_\eps)$.

In view of (\ref{eqtaueps0}), the estimate (\ref{eqM1})
of Lemma~\ref{seclimsup} becomes
\[
\liminf_{\eps\to0}
\eps
\sum_{k=m_\eps}^{n_\eps}
\frac{\PP(\tau_\eps>n_\eps-k)}{A_{k}}
\ge
\frac{\PP(\Omega^*)}{\gamma}
\frac{t}{1+t}.
\]
Combining this with the first part of Lemma~\ref{seclemupperbis}
leads to
\[
\limsup_{\eps\to0} \PP(\tiltau_\eps^x>n_\eps)
\le
\PP(\Gamma_{R,K}^*)
\biggl(
1- \frac{\gamma'}{2f_X(0)} \frac{t}{1+t}
\biggr)
+ o_{R,K}(1).
\]
Successively letting $R\to\infty$, then
$K\to\infty$ and finally $\gamma'\to2f_X(0)$,
we obtain
\[
\limsup_{\eps\to0} \PP(\tiltau_\eps^x>n_\eps)
\le
\frac{\PP(\Omega^*_x)}{1+t}.
\]
To get the corresponding lower bound, parallel to (\ref{eqM2}), we have
\begin{eqnarray*}
&&\PP(\tiltau_\eps^x>2n_\eps) +
{\PP(\Gamma_{R,K}^x) {\gamma''}\eps}
\sum_{k=m_\eps}^{2n_\eps}
\frac{\PP(\tau_\eps>2n_\eps-k)}{A_{k}}\\
&&\qquad\le
\PP(\tiltau_\eps^x>n_\eps) +
t\frac{\gamma''}{\gamma} \PP(\Gamma_{R,K}^x)
{\PP(\tau_\eps>n_\eps)}+o(1).
\end{eqnarray*}
Together with the second part of Lemma~\ref{seclemupperbis}
(with $n_\eps$ replaced by $2n_\eps$)
and (\ref{eqtaueps0}), this implies
\[
\liminf_{\eps\to0} \PP(\tiltau_\eps^x>n_\eps)
\ge
\frac{\PP(\Omega^*_x)}{1+t}
\]
completing the proof.
\end{pf*}

\begin{pf*}{Proof of Theorem~\ref{secdistrib}}
We fix $t \in(0,\infty)$, and choose $n_\eps$
such that $G(n_\eps)\le\frac{\beta t}{\gamma\eps}<G(n_\eps+1)$.

According to the proof of Theorem~\ref{thmtaueps}
[see, in particular, (\ref{eqO}) in Lemma~\ref{seclower}],
we know that for $m_\eps$ with $m_\eps= o(n_\eps)$,
\[
\lim_{\eps\to0}\eps
\sum_{k=m_\eps}^{n_\eps}
\frac{\PP(\tau_\varepsilon>{n_\eps}-k) } {A_{k}}
=
\frac{\PP(\Omega^*)}{\gamma}
\Pr(Y \ge t) =: \psi,
\]
where $Y=\Gamma(\frac{1}{\beta})^{-1} {\mathcal E} {\mathcal
G}_{1/\beta}^{1/\beta}$
is the limiting random variable
of the $\gamma\beta^{-1}\varepsilon G(\tau_\varepsilon)$.
Therefore, Lemma~\ref{seclemupperbis} implies that for
$R,K\ge1$ and $\gamma'<2f_X(0)<\gamma''$,
\[
\limsup_{\varepsilon\rightarrow0}{\mathbb P}(\tiltau_\eps^x>n_\eps)
\le{\mathbb P}(\Gamma_{R,K}^x)(1- \gamma' \psi)
+ o_{R,K}(1)
\]
and
\[
\liminf_{\varepsilon\rightarrow0}{\mathbb P}(\tiltau_\eps^x>n_\eps)
\ge{\mathbb P}(\Gamma_{R,K}^x)
(1- \gamma'' \psi)+o_{R,K}(1).
\]
Since $\lim_{K\rightarrow+\infty}\lim_{R\rightarrow+\infty}
{\mathbb P}(\Gamma_{R,K}^x)={\mathbb P}(\Omega^*_x)$, we get
\[
{\mathbb P}(\Omega^*_x)(1-\gamma''\psi)
\le\liminf_{\varepsilon\rightarrow0}{\mathbb P}(\tiltau_\eps
^x>n_\eps)
\le\limsup_{\varepsilon\rightarrow0}{\mathbb P}(\tiltau_\eps
^x>n_\eps)
\le{\mathbb P}(\Omega^*_x)(1-\gamma'\psi)
\]
and hence
\[
\lim_{\varepsilon\rightarrow0}{\mathbb P}(\tiltau_\eps^x>n_\eps)
= {\mathbb P}(\Omega^*_x)\bigl(1-2f_X(0)\psi\bigr)
= {\mathbb P}(\Omega^*_x) \Pr(Y > t)
\]
as required.
\end{pf*}

\begin{pf*}{Proof of Corollary~\ref{secCOROL}}
This is an $\alpha=2$ case with $A_n=\sqrt{n}$ and
$f_X(0)=\frac{1}{\sqrt{2 \pi}}$. Recalling that
${\mathcal G}_{1/2}=\frac{1}{2 {\mathcal N}^2}$
in distribution
(cf. Example XIII.3.b of~\cite{Feller}) proves our claim.
\end{pf*}

%suskaldyti doi

% imsref loaded by lrinkeviciute, 2011-11-04 15:18:39
%

\printaddresses

\end{document}